\documentclass[11pt]{article}
\usepackage{amsmath}
\usepackage{amsfonts}
\usepackage{amssymb}
\usepackage{amsthm}
\usepackage[margin=1in]{geometry}
\usepackage[alphabetic]{amsrefs}

\newtheorem*{main_HLR_proj}{Theorem 2.3}
\newtheorem*{main_KLR_proj}{Theorem 2.6}
\newtheorem{positivity}{Conjecture}
\newtheorem{main_symmetry}[positivity]{Conjecture}
\newtheorem{main_KtoH}[positivity]{Conjecture}

\newtheorem{HLR_proj}{Theorem}[section]
\newtheorem{KLR_proj}[HLR_proj]{Theorem}
\newtheorem{KtoH_proj}[HLR_proj]{Corollary}
\newtheorem{GRR}[HLR_proj]{Lemma}
\newtheorem{toric_classes}[HLR_proj]{Lemma}
\newtheorem{csm_proj}[HLR_proj]{Proposition}
\newtheorem{mC_proj}[HLR_proj]{Proposition}

\theoremstyle{definition}
\newtheorem{weight}{Definition}[section]
\theoremstyle{plain}
\newtheorem{weight_orth}[weight]{Theorem}
\newtheorem{H_weight_LR}[weight]{Corollary}
\newtheorem{res}[weight]{Lemma}
\newtheorem{csm_proj_weight}[weight]{Proposition}
\newtheorem{proj_like}[weight]{Proposition}
\newtheorem{res_Pieri}[weight]{Lemma}
\newtheorem{csm_Pieri}[weight]{Proposition}
\newtheorem{K_weight_LR}[weight]{Corollary}
\newtheorem{symmetry}{Conjecture}
\newtheorem{KtoH}[symmetry]{Conjecture}

\title{Explicit Formulas for h-Deformed Structure Constants of Grassmannians}
\author{Yiyan Shou}
\date{}

\begin{document}
\maketitle

\begin{abstract}
The Chern-Schwartz-MacPherson (CSM) and motivic Chern (mC) classes of Schubert cells in a Grassmannian are one parameter deformations of the fundamental classes of the Schubert varieties in cohomology and K-theory respectively. Like the fundamental classes, the deformed classes form a basis for the cohomology and K-theory ring of the Grassmannian. The purpose of this paper is to initiate the study of the structure constants associated to the basis CSM and mC classes in terms of the combinatorics of polynomials. First, we prove formulas for the structure constants of projective spaces that involve binomial coefficients. Then, using residue calculus on wieght functions, we describe the structure constants of projective spaces and certain related structure constants of 2-plane Grassmannians as coefficients of explicit polynomials in one variable. Finally, we propose an approach for obtaining more general results in this direction, and make conjectures generalizing the aformentioned results for projective spaces.
\end{abstract}

\section{Introduction}
\subsection{$h$-Deformed Littlewood Richardson Numbers}
Let $R$ be a commutative ring $R$. Given a finite dimensional $R$-algebra $A$ and a $R$-basis $v_1,...,v_n$, we can describe the multiplication operation in $A$ by giving the \emph{structure constants} associated to the basis. Namely, there are unique $c^k_{i,j}\in R$, where $1\leq i,j,k\leq n$, for which $v_i\cdot v_j=\sum_{k=1}^n c^k_{i,j}v_k$. The case where $R=\mathbb{Z}$ and $A$ is the integral cohomology or K-theory ring of a Grassmannian Gr$(d,n)$ has been extensively studied. A Grassmannian is equipped with a canonical CW decomposition into Schubert cells. There are several conventions for indexing these cells. In this paper, Schubert cells $\Omega_I$ will be indexed by subsets $I\subset[1,n]$ such that $|I|=d$. Let $\mathcal{I}=\{I\subset[1,n]\ |\ |I|=d\}$. This immediately gives a basis of fundamental classes of cell closures, $[\overline\Omega_I]$ for $I\in\mathcal{I}$, for cohomology and K-theory. The structure constants associated with this basis are called the \emph{Littlewood-Richardson numbers}. The Littlewood-Richardson numbers are well understood, and several classical rules exist for calculating them. It is possible, however, to consider other bases for the cohomology and K-theory.

The fundamental class admits a 1-parameter deformation. In cohomology, the deformation is called the \emph{Chern-Schwartz-MacPherson} or $c^{sm}$ class \cite{M}. In K-theory, it is called the \emph{motivic Chern} or $mC$ class \cite{BSY}. Throughout this paper, this parameter will be called $h$. If $X$ is a smooth variety, these classes associate an element of $\mathrm{H}^*(X)[h]$ and $\mathrm{K}^0(X)[h]$ respectively to each subvariety of $X$. Taking the term of highest $h$ power in the $\mathrm{c^{sm}}$ class of a locally closed subvariety recovers the cohomological fundamental class of the closure, and setting $h$ to zero in mC recovers the K-theoretic fundamental class. Alternatively, one can think of $\mathrm{c^{sm}}$ and mC as generalizing the total Chern class of the tangent bundle of a smooth variety to potentially singular varieties. In particular, integrating these classes gives the Euler characteristic (up to a power of $h$) and $\chi_h$ characteristic respectively. An important property of these classes is \emph{additivity}: given disjoint subvarieties $S,T$,
	\[\mathrm{c^{sm}}(S\cup T)=\mathrm{c^{sm}}(S)+\mathrm{c^{sm}}(T)\ \ ,\ \ \mathrm{mC}(S\cup T)=\mathrm{mC}(S)+\mathrm{mC}(T).\]
This paper deals with the $\mathrm{c^{sm}}$ and mC classes of Schubert cells.

The $h$-deformed classes of Schubert cells in partial flag varieties are studied in \cite{AMSS1,AMSS2}. Like the fundamental classes, the $h$-deformed classes of the Schubert cells form a basis for cohomology and K-theory. An immediate consequence of the above properties is
	\[\forall I\in\mathcal{I}\ \ \mathrm{c^{sm}}(\Omega_I)=h^{\dim(\Omega_I)}[\overline\Omega_I]+\text{h.o.t},\]
where "h.o.t" denotes terms of higher cohomological degree. We will refer to this property as \emph{triangularity with respect to cohomological grading}. By performing a procedure akin to Gaussian elimination over $\mathbb{Z}(h)$, we can uniquely express the fundamental classes $[\overline\Omega_I]$ in terms of the Chern-Schwartz-MacPherson classes $\mathrm{c^{sm}}(\Omega_I)$. Therefore, the classes $\mathrm{c^{sm}}(\Omega_I)$, where $I\in\mathcal{I}$, form a $\mathbb{Z}(h)$-basis for $\mathrm{H}^*(\mathrm{Gr}(d,n))(h)$. A similar argument with Chern character shows that the classes $\mathrm{mC}(\Omega_I)$ for $I\in\mathcal{I}$ form a $\mathbb{Z}(h)$-basis for $\mathrm{K}^0(\mathrm{Gr}(d,n))(h)$. The associated structure constants $c^K_{I,J}$ and $C^K_{I,J}$ in cohomology and K-theory respectively are called the \emph{h-deformed Littlewood-Richardson numbers}. A formula for the cohomological h-deformed structure constants of a partial flag variety in terms of divided difference operators is given in \cite{S}. It is also possible to consider the structure constants associated to the bases of Segre-Schwartz-MacPherson and Segre-Motivic Chern classes in cohomology and K-theory respectively. These structure constants are closely related to our structure constants (see Section 2.4 of \cite{S}).

In this paper, we calculate the $h$-deformed Littlewood-Richardson numbers of projective spaces, give preliminary results for 2-plane Grassmannians, and propose an approach for generalizing these results to arbitrary Grassmannians. The main results and conjectures are listed below. In Section 2, we will consider the $h$-deformed Littlewood Richardson numbers of projective spaces $\mathrm{Gr}(1,m+1)=\mathbb{P}^m$. Our approach is based on toric calculus, and in Section 3, we formulate the problem of finding $h$-deformed structure constants for arbitrary toric manifolds. Since the Schubert cells of $\mathbb{P}^m$ are indexed by singleton sets, we drop the set braces from the notation. Hence, the cells will be indexed by integers $1\leq i\leq m+1$. In Section 4, we give some preliminary results and conjectures for Grassmannians that are not projective spaces. Our approach is based on using weight function orthogonality \cite{R,RTV1,RTV2} and residue calculus (cf. \cite{FR}) to obtain formulas for the equivariant $h$-deformed Littlewood-Richardson numbers and then degenerating to the nonequivariant structure constants by taking limits.

\subsection*{Acknowledgements}
I would like to thank my advisor, Rich\'ard Rim\'anyi, for his invaluable guidance and comments on this paper.

\subsection{Results and Conjectures}
\begin{main_HLR_proj}
Let $c^k_{i,j}$ be the structure constants of $\mathbb{P}^m$ associated to the basis of CSM classes of Shubert cells. For all $1\leq i,j,k\leq m+1$,
\begin{enumerate}
\item $c^k_{i,j}=0$ unless $k\leq i,j$,
\item $c^k_{i,j}=c^k_{i-1,j+1}$ when $i>1$ and $j<m+1$,
\item $c^k_{i,j}=c^{k+1}_{i,j+1}$ when $j,k<m+1$,
\item $c^k_{m+1,m+1}=\binom{2m-k}{m-1}h^m$.
\end{enumerate}
These equalities completely determine all cohomological structure constants of $\mathbb{P}^m$. In particular,
	\[c^k_{i,j}=\binom{i+j-k-2}{m-1}h^m.\]
\end{main_HLR_proj}

\begin{main_KLR_proj}\label{KLR_proj}
Let $C^k_{i,j}$ be the structure constants of $\mathbb{P}^m$ associated to the basis of motivic Chern classes of Schubert cells. Define $\widetilde{C}^K_{i,j}=(-1)^{m+i+j+k+1}(1+h)^{-m}C^K_{i,j}$. For all $1\leq i,j,k\leq m+1$,
\begin{enumerate}
\item $C^k_{i,j}=0$ unless $k\leq i,j$,
\item $C^k_{i,j}=C^k_{i-1,j+1}$ when $i>1$ and $j<m+1$,
\item $C^k_{i,j}=C^{k+1}_{i,j+1}$ when $j,k<m+1$,
\item $\widetilde{C}^k_{m+1,m+1}=\binom{2m-k}{m}h^{m-k}+\binom{2m-k+1}{m}h^{m-k+1}$.
\end{enumerate}
These equalities completely determine all K-theoretic structure constants of $\mathbb{P}^m$. In particular,
	\[\widetilde C^k_{i,j}=\binom{i+j-k-2}{m}h^{i+j-k-m-2}+\binom{i+j-k-1}{m}h^{i+j-k-m-1}.\]
\end{main_KLR_proj}

\begin{main_symmetry}
For all $I,I',J,J',K,K'\in\mathcal{I}$,
\begin{enumerate}
\item $C^K_{I,J}$ is a polynomial in $-1-h$ with nonnegative coefficients,
\item $C^K_{I,J}=0$ unless $K\leq I,J$,\\
\item $C^K_{I,J}=C^K_{I',J'}$ if $i_a+j_a=i'_a+j'_a$ for all $1\leq a\leq d$,\\
\item $C^K_{I,J}=C^{K'}_{I,J'}$ if $j_a-k_a=j'_a-k'_a$ for all $1\leq a \leq d$.
\end{enumerate}
\end{main_symmetry}

\begin{main_KtoH}
Expand $C^K_{I,J}$ as a polynomial in $-1-h$. Then the lowest degree term has coefficient $c^K_{I,J}/h^{\dim(\mathrm{Gr}(d,n))}$.
\end{main_KtoH}



\section{Projective Spaces}
\subsection{Toric Manifolds}
In Section 2, we will compute structure constants using toric geometry. Our reference for toric geometry is \cite{Fu}. First, let us recall the necessary constructions. The geometry of a toric manifold is encoded in a Delzant polytope. A convex polytope $\Delta\subset\mathbb{R}^m$ is Delzant if
	\begin{enumerate}
	\item there are $m$ edges incident to each vertex,
	\item each edge incident to a vertex $p$ is of the form $\{p+tu\ |\ t\in[0,h]\}$, where $h\in\mathbb{R},u\in\mathbb{Z}^m\subset\mathbb{R}^m$,
	\item and these integral vectors $u$ form a $\mathbb{Z}$-basis for $\mathbb{Z}^m$.
	\end{enumerate}
Let $T=(\mathbb{C}^*)^m$ be a torus and $X$ be a smooth compact complex Hamiltonian $T$-space. This action restricted to $U=U(1)^m$ is Hamiltonian, and we denote the moment map by $\mu$. We say $X$ is a toric manifold if the action of $U$ is effective and $m=\mathrm{dim}X$ (all dimensions in this paper are complex dimensions). In this case, $\Delta=\mu(X)\subset\mathrm{Lie}(U)=\mathbb{R}^m$ is a Delzant polytope, and there is a correspondence between $d$-faces of $\Delta$ and $d$-dimensional orbits of $T$. We will take the dual perspective by looking at the normal fan to $\Delta$. This forgets the edge lengths and position of $\Delta$ and records only the angles between facets, but the forgotten data encodes the symplectic structure, which is not relevant for this paper.

Let $F_1,...,F_r$ be the facets of $\Delta$. Let $v_i\in\mathbb{Z}^m$ be the primitive inward pointing normal vector to $F_i$. For $S\subset [1,r]$, define the cone $\sigma_S=\sum_{i\in S}\mathbb{R}_{\geq 0}\cdot v_i$. The normal fan to $\Delta$ is the collection of cones $\mathcal{N}(\Delta)=\{\sigma_S|\bigcap_{i\in S}F_i\neq\emptyset\}$. The zero cone $\{0\}$ belongs to $\mathcal{N}(\Delta)$ by convention. In the language of fans, the Delzant property translates to
	\begin{enumerate}
	\item the maximal cones are spanned by $m$ vectors,
	\item each cone is spanned by integral vectors,
	\item and the vectors spanning each maximal cone form a $\mathbb{Z}$-basis for $\mathbb{Z}^m$.
	\end{enumerate}
Moreover, there is a correspondence between $d$-dimensional faces of $\Delta$ and $(m-d)$-dimensional cones of $\mathcal{N}(\Delta)$. This gives a correspondence of $d$-dimensional cones and $T$-orbits of codimension $d$
	\[\{\gamma\in\mathcal{N}(\Delta)\ |\ \dim(\gamma)=d\} \longleftrightarrow \{T\text{-orbits of codimension }d\}.\]
Under this correspondence, fixed points are identified with maximal cones, 1-dimensional orbits are identified with facets, and codimension 1 orbits are identified with rays. The orbit structure of $X$ has a simple combinatorial description.

Given $\gamma\in\mathcal{N}(\Delta)$, let $\mathcal{O}_\gamma$ be the orbit corresponding to $\gamma$, and $V_\gamma=\overline{\mathcal{O}_\gamma}$. Then, $V_\gamma=\bigcup_{\gamma'\supset\gamma}\mathcal{O}_{\gamma'}$. It follows that the correspondence between cones and orbits is inclusion reversing with respect to orbit closures. It also follows that if $\gamma_1,\gamma_2,\gamma_3\in\mathcal{N}(\Delta)$ are such that $\gamma_3=\gamma_1+\gamma_2$, then $V_{\gamma_3}=V_{\gamma_1}\cap V_{\gamma_2}$. If in addition $\gamma_1\cap\gamma_2=0$, then a local coordinates calculation shows that this intersection is transverse. In the case that $\gamma_1$ and $\gamma_2$ do not span a cone in $\mathcal{N}(\Delta)$, $V_{\gamma_1}\cap V_{\gamma_2}=\emptyset$. All orbit closures $V\gamma$ are themselves toric manifolds, and their normal cones are given by the star construction. Define the star of $\gamma$, to be the following fan of cones in $\mathbb{R}^m/(\mathbb{R}\cdot\gamma)\cong\mathbb{R}^{m-\dim(\gamma)}$:
	\[\mathrm{st}(\gamma):=\{\gamma'+\mathbb{R}\cdot\gamma\in\mathbb{R}^m/(\mathbb{R}\cdot\gamma)\ |\ \gamma'\in\mathcal{N}(\Delta),\gamma'\supset\gamma\}.\]
The moment polytope of $V_\gamma$ is the corresponding face of $\Delta$. These properties give rise to a combinatorial description of $\mathrm{H}^*(X)$.

Let $\tau_1,...,\tau_r$ be the rays in $\mathcal{N}(\Delta)$. Define the Stanley-Reisner ring of $\Delta$ by 
	\[SR:=\mathbb{Z}[X_1,...,X_r]/(X_{i_1}\cdots X_{i_s}\ |\ s\leq r,1\leq i_a\leq r,\text{ and }\tau_{i_1},...,\tau_{i_s}\text{ do not span a cone in }\mathcal{N}(\Delta)).\]
We will see that $\mathrm{H}^*(X)\cong SR/\mathcal{J}$, where $\mathcal{J}=\left\langle\sum_{i=1}^r u(\tau_i)X_i\ \middle|\ {u\in\mathrm{Hom}(\mathbb{Z}^m,\mathbb{Z})}\right\rangle$. By assembling torus orbits into cells, we can build a CW decomposition of $X$ (more on this construction in Section 3), where the closure of each cell is a torus orbit closure $V_\gamma$ for some $\gamma\in\mathcal{N}(\Delta)$. This CW decomposition implies that $\mathrm{H}^*(X)$ is generated by the fundamental classes $[V_\gamma]$. If $\gamma$ is spanned by distinct rays $\tau_{i_1},...,\tau_{i_s}$, then the transverse intersection property implies that $[V_\gamma]=[V_{\tau_{i_1}}]\cdots[V_{\tau_{i_s}}]$. Moreover, if $\tau_{i_1},...,\tau_{i_s}$ do not span a cone in $\mathcal{N}(\Delta)$, then $V_{\tau_{i_1}}\cap\cdots\cap V_{\tau_{i_s}}=\emptyset$, and $[V_{\tau_{i_1}}]\cdots[V_{\tau_{i_s}}]=0$. Hence, there is a surjection
\begin{align*}
\Phi:SR&\to \mathrm{H}^*(X)\\
X_i&\mapsto[V_{\tau_i}]
\end{align*}
The divisors $\sum_{i=1}^r u(\tau_i)V_{\tau_i}$ are known to be principal Cartier, so certainly $\mathcal{J}\subset\ker(\Phi)$. It is shown in \cite{Fu1}, that in fact $\mathcal{J}=\ker(\Phi)$, giving the desired isomorphism $SR/\mathcal{J}\cong \mathrm{H}^*{X}$. We will be interested in the CSM and motivic Chern classes of the cells in the CW decomposition.

In order to compute motivic characteristic classes of a cell, we can stratify it into $T$-orbits and apply the additivity property. The following formulas for the CSM and motivic Chern class of a $T$-orbit are given in \cite{MS}. For all $\gamma\in\mathcal{N}(\Delta)$,
	\[\mathrm{c^{sm}}(\mathcal{O}_\gamma)=h^{m-\dim(\gamma)}[V_\gamma],\quad\mathrm{mC}_h(\mathcal{O}_\gamma)=(1+h)^{m-\dim(\gamma)}\iota_{\gamma*}\omega_\gamma,\]
where $\omega_\gamma$ is the canonical bundle of $V_\gamma$, and $\iota_\gamma:V_\gamma\to X$ is the inclusion. It is not necessary to describe the K-theory of $X$, as we will immediately pass to cohomology by taking Chern characters. The following lemma records some elementary facts that are needed for future calculations.

\begin{toric_classes}\label{toric_classes}
Fix $\gamma\in\mathcal{N}(\Delta)$. Given $\eta\in\mathcal{N}(\Delta)$ such that $\eta\supset\gamma$, let $\tilde\eta=\eta+\mathbb{R}\cdot\gamma\in\mathrm{st}(\gamma)$, we have
\begin{enumerate}
\item $\begin{displaystyle} V_{\{0\}}=X \end{displaystyle}$,
\item $\begin{displaystyle} \mathrm{c}^*(T_*V_\gamma)=\prod_{\dim(\tau)=2}(1+[V_{\tilde\tau}]) \end{displaystyle}$,
\item $\begin{displaystyle} \mathrm{c}^*(\omega_\gamma)=-\sum_{\dim(\tau)=2}[V_{\tilde\tau}] \end{displaystyle}$,
\item $\begin{displaystyle} \iota_{\gamma*}[V_{\tilde\eta}]=[V_\eta] \end{displaystyle}$.
\end{enumerate}
\end{toric_classes}

The projective space $\mathbb{P}^m$ has an especially simple toric description. Note that the usual torus action on $\mathbb{P}^m$, induced by the action on $\mathbb{C}^{m+1}$, is not effective. However, we resolve this by quotienting by the kernel of the action. Equivalently, we take $T$ to be the first $m$ factors of the usual $(m+1)$-dimensional torus, so that $T$ acts trivially on the $(m+1)$th homogeneous coordinate. Under standard conventions for the moment map, the moment polytope $\Delta$ is the convex hull of $0,(1/2)e_1,...,(1/2)e_m$, where the $e_i$ are the standard basis vectors of $\mathbb{R}^m$. Defining $e_{m+1}=-e_1-\cdots-e_m$, we have
	\[\mathcal{N}(\Delta)=\{\mathbb{R}_{\geq 0}\cdot S|S\subsetneq\{e_1,...,e_{m+1}\}\}\cup\{\{0\}\}.\]
Let $\sigma_i=\mathbb{R}_{\geq 0}\cdot\{e_1,...,\hat{e_i},...,e_{m+1}\}$ for $1\leq i\leq m+1$. These are the maximal cones, which correspond to the $m+1$ fixed points of $\mathbb{P}^m$. Let $\tau_i=\mathbb{R}_{\geq 0}\cdot e_i$ for $1\leq i\leq m+1$. These are the rays, which correspond to toric divisors. It is easy to see that the divisors $V_{\tau_i}$ are all equivalent under the relation of $\mathcal{J}$. Therefore, we denote the common value of $[V_{\tau_i}]$ for $1\leq i\leq m+1$, by $H$. In fact, $H\in\mathrm{H}^2(\mathbb{P}^m)$ is the hyperplane class.
Define $\gamma_S=\mathbb{R}_{\geq 0}\cdot \{e_i\ |\ i\in S\}\in\mathcal{N}(\Delta)$, for $S\subsetneq[1,m+1]$. For $1\leq i\leq m+1$, define $\Omega_i\subset\mathbb{P}^m$ inductively as follows:
\begin{align*}
\Omega_1&=\mathcal{O}_{\gamma_{[1,m]}},\\
\Omega_i&=\bigcup_{S\supset[1,m-i+1]}\mathcal{O}_{\gamma_S}\setminus\Omega_{i-1}\text{, when }i>1.\\
\end{align*}
Then, $\Omega_i\cong\mathbb{C}^{i-1}$, and the $\Omega_i$ decompose $\mathbb{P}^m$ into Schubert cells. The number of $d$-dimensional orbits in the stratification of $\Omega_i$ is $\binom{i-1}{d}$.

\subsection{Cohomological Structure Constants of $\mathbb{P}^m$}
The classes $\mathrm{c^{sm}}(\Omega_i)$ for $1\leq i\leq m+1$ form a basis for $\mathrm{H}^*(\mathbb{P}^m)$. Let $c^k_{i,j}\in\mathbb{Z}(h)$ be the unique rational functions satisfying $\forall (1\leq i,j\leq m+1)\ \ \mathrm{c^{sm}}(\Omega_i)\mathrm{c^{sm}}(\Omega_j)=\sum_{k=1}^{m+1}c^k_{i,j}\mathrm{c^{sm}}(\Omega_k)$. Given $\gamma\in\mathcal{N}(\Delta)$, our description of the cohomology of $\mathbb{P}^m$ gives the formula
	\[\mathrm{c^{sm}}(\mathcal{O}_\gamma)=h^{m-\dim(\gamma)}[V_\gamma]=h^{m-\dim(\gamma)}\cdot\prod_{\tau\text{ an extremal ray of }\gamma}[V_\tau]=h^{m-\dim(\gamma)}[H]^{\mathrm{dim}(\gamma)}.\]

\begin{csm_proj}\label{csm_proj}
For all $1\leq i\leq m+1$,\ \ $\mathrm{c^{sm}}(\Omega_i)=H^{m-i+1}(h+H)^{i-1}$.
\end{csm_proj}

\begin{proof}
The additivity of $\mathrm{c^{sm}}$ classes implies 
\begin{align*}
\mathrm{c^{sm}}(\Omega_i)&=\sum_{\substack{\gamma\in\mathcal{N}(\Delta)\\ \mathcal{O}_\gamma\subset\Omega_i}}\mathrm{c^{sm}}(\mathcal{O}_\gamma)\\
	&=\sum_{\substack{\gamma\in\mathcal{N}(\Delta)\\ \mathcal{O}_\gamma\subset\Omega_i}}h^{m-\dim(\gamma)}H^{\mathrm{dim}(\gamma)}.
\end{align*}
Because $\Omega_i$ contains $\binom{i-1}{d}$ many orbits of dimension $d$, the sum is over $\binom{i-1}{d}$ cones of dimension $m-d$. All cones in the sum have dimension at least $m-i+1$. We have,
\begin{align*}
\mathrm{c^{sm}}(\Omega_i)&=\sum_{d=m-i+1}^{m}\binom{i-1}{m-d}h^{m-d}H^d\\
	&=H^{m-i+1}\sum_{d=0}^{i-1}\binom{i-1}{d}h^{i-1-d}H^d\\
	&=H^{m-i+1}(h+H)^{i-1}.
\end{align*}
\end{proof}

\begin{HLR_proj}\label{HLR_proj}
For all $1\leq i,j,k\leq m+1$,
\begin{enumerate}
\item $c^k_{i,j}=0$ unless $k\leq i,j$,
\item $c^k_{i,j}=c^k_{i-1,j+1}$ when $i>1$ and $j<m+1$,
\item $c^k_{i,j}=c^{k+1}_{i,j+1}$ when $j,k<m+1$,
\item $c^k_{m+1,m+1}=\binom{2m-k}{m-1}h^m$.
\end{enumerate}
These equalities completely determine all cohomological structure constants of $\mathbb{P}^m$. In particular,
	\[c^k_{i,j}=\binom{i+j-k-2}{m-1}h^m.\]
\end{HLR_proj}

\begin{proof}
\begin{enumerate}
\item This follows from degree considerations.

\item  Using Proposition~\ref{csm_proj}, we have
\begin{align*}
\mathrm{c^{sm}}(\Omega_i)\mathrm{c^{sm}}(\Omega_j)&=H^{m-i+1}(h+H)^{i-1}\cdot H^{m-j+1}(h+H)^{j-1}\\
	&=H^{2m-i-j+2}(h+H)^{i+j-2}\\
	&=H^{m-(i-1)+1+m-(j+1)+1}(h+H)^{(i-1)-1+(j+1)-1}\\
	&=H^{m-(i-1)+1}(h+H)^{(i-1)-1}\cdot H^{m-(j-1)+1}(h+H)^{(j-1)-1}\\
	&=\mathrm{c^{sm}}(\Omega_{i-1})\mathrm{c^{sm}}(\Omega_{j+1}).
\end{align*}

\item We have $\mathrm{c^{sm}}(\Omega_{i})\mathrm{c^{sm}}(\Omega_{j})=H^{2m-i-j+2}(h+H)^{i+j-2}$, and $\mathrm{c^{sm}}(\Omega_{i})\mathrm{c^{sm}}(\Omega_{j+1})=H^{2m-i-j+1}(h+H)^{i+j-1}$ from Proposition 2.2. By part (1), the $c^k_{i,j}$ satisfy
\begin{align*}
 H^{2m-i-j+2}(h+H)^{i+j-2}&=\sum_{k=1}^{j}c^k_{i,j}\mathrm{c^{sm}}(\Omega_k)\\
	&=\sum_{k=1}^{j}c^k_{i,j}H^{m-k+1}(h+H)^{k-1}.
\end{align*}
Multiplying both sides by $(h+H)$ yields
	\[H\cdot H^{2m-i-j+1}(h+H)^{i+j-1}=\sum_{k=1}^{j}c^k_{i,j}H\cdot H^{m-k}(h+H)^{k}.\]
Applying Proposition~\ref{csm_proj}, we get
	\[H\cdot\mathrm{c^{sm}}(\Omega_{i})\mathrm{c^{sm}}(\Omega_{j+1})=H\cdot\sum_{k=1}^{j}c^k_{i,j}\mathrm{c^{sm}}(\Omega_{k+1}).\]
The factor of $H$ kills the top degree part, but we have equality in all other degrees. Because the basis of $\mathrm{c^{sm}}$ classes is triangular with respect to cohomological grading, the top degree part is controlled by the $c^1_{i,j+1}$ structure constant.

\item We will calculate $\sum_{k=1}^{m+1}\binom{2m-k}{m-1}h^m\mathrm{c^{sm}}(\Omega_k)$, and show that it is equal to $\mathrm{c^{sm}}(\Omega_{m+1})\mathrm{c^{sm}}(\Omega_{m+1})$. We have
\begin{align*}
\sum_{k=1}^{m+1}\binom{2m-k}{m-1}h^m\mathrm{c^{sm}}(\Omega_k)&=\sum_{k=1}^{m+1}\binom{2m-k}{m-1}h^m\cdot H^{m-k+1}(h+H)^{k-1}\\
	&=\sum_{k=1}^{m+1}\binom{2m-k}{m-1}h^mH^{m-k+1}\sum_{l=0}^{k-1}h^lH^{k-l-1}\binom{k}{l}\\
	&=\sum_{k=1}^{m+1}\sum_{l=0}^{k-1}h^{m+l}H^{m-l}\binom{2m-k}{m-1}\binom{k-1}{l}\\
	&=\sum_{l=0}^{m}h^{m+l}H^{m-l}\sum_{k=l+1}^{m+1}\binom{2m-k}{m-1}\binom{k-1}{l}\\
	&=\sum_{l=0}^{m}h^{m+l}H^{m-l}\sum_{k=0}^{m-l}\binom{2m-(k+l)-1}{m-1}\binom{k+l}{l}.
\end{align*}
On the other hand, because $H^l=0$ when $l>m$,
\begin{align*}
\mathrm{c^{sm}}(\Omega_{m+1})\mathrm{c^{sm}}(\Omega_{m+1})&=(h+H)^{2m}\\
	&=\sum_{l=0}^{2m}h^{2m-l}H^{l}\binom{2m}{l}\\
	&=\sum_{l=0}^{m}h^{2m-l}H^{l}\binom{2m}{l}\\
	&=\sum_{l=0}^{m}h^{m+l}H^{m-l}\binom{2m}{m-l}\\
	&=\sum_{l=0}^{m}h^{m+l}H^{m-l}\binom{2m}{m+l}.
\end{align*}
By comparing the powers of $h$, it suffices to show $\binom{2m}{m+l}=\sum_{k=0}^{m-l}\binom{2m-(k+l)-1}{m-1}\binom{k+l}{l}$ for all $0\leq l\leq m$. This is precisely the Vandermonde identity.
\end{enumerate}
\end{proof}

\subsection{K-Theoretic Structure Constants of $\mathbb{P}^m$}
As before, the classes $\mathrm{mC}(\Omega_i)$ for $1\leq i\leq m+1$ form a basis for $K^0(\mathbb{P}^m)$. Let $C^k_{i,j}\in\mathbb{Z}(h)$ be the unique rational functions satisfying $\forall (1\leq i,j\leq m+1)\ \ \mathrm{mC}(\Omega_i)\mathrm{mC}(\Omega_j)=\sum_{k=1}^{m+1}C^k_{i,j}\mathrm{mC}(\Omega_k)$. We will not describe the K-theory ring of $\mathbb{P}^m$. Instead, we use the fact that any K-theory class is uniquely determined by its Chern character. Our main tool will be the following corollary of the Grothendieck-Riemann-Roch theorem.

\begin{GRR}\label{GRR}
Let $X$ be a smooth projective variety and $V\subset X$ a smooth closed subvariety. Denote the inclusion by $\iota:V\to X$. Let $\xi\to V$ be a vector bundle. Then,
	\[\mathrm{ch}(\iota_*\xi)=\iota_*(\mathrm{ch}(\xi)\mathrm{td}(V))/\mathrm{td}(X).\]
\end{GRR}

Fix $\gamma\in\mathcal{N}(\Delta)$. Given $\eta\in\mathcal{N}(\Delta)$ with $\eta\supset\gamma$, let $\tilde\eta=\eta/\gamma\in\mathrm{st}(\gamma)$. We begin by computing the Chern character of pushforwards of the canonical bundle $\omega_\gamma\to V_\gamma$. By Lemma 2.1 (4), we have $\forall\eta\in\mathcal{N}(\Delta)\ \ [V_{\tilde\eta}]=H_\gamma^{\dim(\tilde\eta)}\in \mathrm{H}^*(V_\gamma)$, where $H_\gamma$ is the hyperplane class of $V_\gamma$. This along with parts (2) and (3) of Lemma~\ref{toric_classes} yield the formulas
	\[\mathrm{ch}(\omega_\gamma)=\mathrm{e}^{-(m-\dim(\gamma)+1)H_\gamma},\quad\mathrm{td}(V_\gamma)=\left(\frac{H_\gamma}{1-\mathrm{e}^{-H_\gamma}}\right)^{m-\dim(\gamma)+1}.\]
Using Lemma~\ref{GRR} and part (4) of Lemma~\ref{toric_classes}, we have
\begin{align*}
\mathrm{ch}(\iota_{\gamma*}\omega_\gamma)&=\iota_{\gamma*}(\mathrm{ch}(\omega_\gamma)\mathrm{td}(V_\gamma))/\mathrm{td}(\mathbb{P}^m)\\
	&=H^{\dim(\gamma)}\cdot\frac{\mathrm{e}^{-(m-\dim(\gamma)+1)H}H^{m-{\dim(\gamma)}+1}}{(1-\mathrm{e}^{-H})^{m-\dim(\gamma)+1}}\Big/\frac{H^{m+1}}{(1-\mathrm{e}^{-H})^{m+1}}\\
	&=\mathrm{e}^{-(m-\dim(\gamma)+1)H}(1-\mathrm{e}^{-H})^{\dim(\gamma)}.
\end{align*}
Combining this formula with the orbit structure of $\Omega_i$, the additivity property of motivic Chern classes, and naturality of Chern character results in the following formula.

\begin{mC_proj}\label{mC_proj}
For all $1\leq i\leq m+1$,\ \ $\mathrm{ch}(\mathrm{mC}(\Omega_i))=\mathrm{e}^{-H}(1-\mathrm{e}^{-H})^{m-i+1}(1+h\mathrm{e}^{-H})^{i-1}$.
\end{mC_proj}

The proof is similar to that of Proposition~\ref{csm_proj}. Notice the similarities between the formulas of Proposition~\ref{mC_proj} and Proposition~\ref{csm_proj}.

\begin{KLR_proj}\label{KLR_proj}
Let $\widetilde{C}^k_{i,j}=(-1)^{m+i+j+k+1}(1+h)^{-m}C^k_{i,j}$. For all $1\leq i,j,k\leq m+1$,
\begin{enumerate}
\item $C^k_{i,j}=0$ unless $k\leq i,j$,
\item $C^k_{i,j}=C^k_{i-1,j+1}$ when $i>1$ and $j<m+1$,
\item $C^k_{i,j}=C^{k+1}_{i,j+1}$ when $j,k<m+1$,
\item $\widetilde{C}^k_{m+1,m+1}=\binom{2m-k}{m}h^{m-k}+\binom{2m-k+1}{m}h^{m-k+1}$.
\end{enumerate}
These equalities completely determine all K-theoretic structure constants of $\mathbb{P}^m$. In particular,
	\[\widetilde C^k_{i,j}=\binom{i+j-k-2}{m}h^{i+j-k-m-2}+\binom{i+j-k-1}{m}h^{i+j-k-m-1}.\]
\end{KLR_proj}

\begin{proof}
The first three parts can be proven the same way the first three parts of Theorem~\ref{HLR_proj} were. It remains to prove part 4. Observe that $\rho:=1-\mathrm{e}^{-H}=H-1/2H^2+\mathrm{h.o.t}$. It follows that for $0\leq l\leq m$, the elements $\rho^l$ are triangular with respect to the grading on $\mathrm{H}^*(\mathbb{P}^m)$. Thus, these elements form a basis for the cohomology ring. We can express Proposition~\ref{mC_proj} in terms of $\rho$ as
	\[\mathrm{ch}(\mathrm{mC}(\Omega_i))=(1-\rho)\rho^{m-i+1}((1+h)-h\rho)^{i-1}\quad (*).\]
We will expand $h^{-m}\mathrm{ch}(\mathrm{mC}(\Omega_{m+1}))^2$ and $\sum_{k=1}^{m}(-1)^{m+k+1}\left(\binom{2m-k}{m}h^{m-k}+\binom{2m-k+1}{m}h^{m-k+1}\right)\mathrm{ch}(\mathrm{mC}(\Omega_k))$ in powers of $\rho$, using $(*)$, and compare coefficients. We have
\begin{align*}
h^{-m}\mathrm{ch}(\mathrm{mC}(\Omega_{m+1}))^2&=h^{-m}(1-\rho)^2((1+h)-h\rho)^{2m}\\
&=\sum_{l=0}^m\rho^l(-1)^l\left(\binom{2m}{l}h^{l}(1+h)^{m-l}+2\binom{2m}{l-1}h^{l-1}(1+h)^{m-l+1}\right.\\
&\quad\qquad\qquad\qquad\left.+\binom{2m}{l-2}h^{l-2}(1+h)^{m-l+2}\right).
\end{align*}
On the other hand,
\begin{align*}
\sum_{k=1}^{m}(-1&)^{m+k+1}\left(\binom{2m-k}{m}h^{m-k}+\binom{2m-k+1}{m}h^{m-k+1}\right)\mathrm{ch}(\mathrm{mC}(\Omega_k))\\
	&=\sum_{k=1}^{m}(-1)^{m+k+1}\left(\binom{2m-k}{m}h^{m-k}+\binom{2m-k+1}{m}h^{m-k+1}\right)\\
		&\hspace{30pt}\cdot(1-\rho)\rho^{m-k+1}((1+h)-h\rho)^{k-1}\\
	&=\sum_{l=0}^m\rho_l\sum_{k=m-l+1}^{m+1}(-1)^{m+k+1}\left(\binom{2m-k}{m}h^{m-k}+\binom{2m-k+1}{m}h^{m-k+1}\right)\\
		&\hspace{70pt}\cdot\left(\binom{k-1}{l+k-m-1}(-h)^{l+k-m-1}(1+h)^{m-l}\right.\\
			&\hspace{96pt}\left.-\binom{k-1}{l+k-m-2}(-h)^{l+k-m-2}(1+h)^{m-l+1}\right)\\
	&=\sum_{l=0}^m\rho_l(-1)^l\sum_{k=m-l+1}^{m+1}\left(\binom{2m-k}{m}\binom{k-1}{m-l}h^{l-1}(1+h)^{m-l}\right.\\
		&\hspace{96pt}+\binom{2m-k}{m}\binom{k-1}{m-l+1}h^{l-2}(1+h)^{m-l+1}\\
		&\hspace{96pt}+\binom{2m-k+1}{m}\binom{k-1}{l-m}h^l(1+h)^{m-l}\\
		&\left.\hspace{96pt}+\binom{2m-k+1}{m}\binom{k-1}{m-l+1}h^{l-1}(1+h)^{m-l+1}\right).
\end{align*}
Consider the final sum. Apply the Vandermonde identity to each of the four lines, and then apply Pascal's identity only to the last two lines. This yields
\begin{align*}
\sum_{k=m-l+1}^{m+1}\binom{2m-k}{m}\binom{k-1}{m-l}h^{l-1}(1+h)^{m-l}&=\binom{2m}{l-1}h^{l-1}(1+h)^{m-l},
\end{align*}
\begin{align*}
\sum_{k=m-l+1}^{m+1}\binom{2m-k}{m}\binom{k-1}{m-l+1}h^{l-2}(1+h)^{m-l+1}&=\binom{2m}{l-2}h^{l-2}(1+h)^{m-l+1},
\end{align*}
\begin{align*}
\sum_{k=m-l+1}^{m+1}\binom{2m-k+1}{m}\binom{k-1}{m-l}h^l(1+h)^{m-l}&=\binom{2m+1}{l}h^l(1+h)^{m-l}\\
	&=\left(\binom{2m}{l}+\binom{2m}{l-1}\right)h^l(1+h)^{m-l},
\end{align*}
\begin{align*}
\sum_{k=m-l+1}^{m+1}\binom{2m-k+1}{m}\binom{k-1}{m-l+1}h^{l-1}(1+h)^{m-l+1}&=\binom{2m+1}{l-1}h^{l-1}(1+h)^{m-l+1}\\
	&=\left(\binom{2m}{l-1}+\binom{2m}{l-2}\right)h^{l-1}(1+h)^{m-l+1}.
\end{align*}
Adding these four expressions yields 
	\[\binom{2m}{l}h^{l}(1+h)^{m-l}+2\binom{2m}{l-1}h^{l-1}(1+h)^{m-l+1}+\binom{2m}{l-2}h^{l-2}(1+h)^{m-l+2},\]
as required.
\end{proof}

The following corollary is a consequence of Pascal's identity, Theorem~\ref{HLR_proj}, and Theorem~\ref{KLR_proj}.

\begin{KtoH_proj}\label{KtoH_proj}
Expand $C^k_{i,j}$ as a polynomial in $-1-h$. Then, the coefficient of the lowest degree term is $c^k_{i,j}/h^m$.
\end{KtoH_proj}

\section{Toric Calculus}
In this section, we will formulate the problem of finding structure constants for arbitrary toric manifolds $X$ with moment polytope $\Delta$. We will then use the machinery developed in Section 2 in the context of projective spaces to calculate the structure constants of other toric manifolds. The first necessary ingredient is a CW decomposition of $X$.

\subsection{Shellings and CW Decompositions}
Let $x_1,...,x_r\in\mathrm{Lie}(U)^*\cong\mathbb{R}^n$ be the vertices of $\Delta$. Fix a vector $u\in\mathrm{Lie}(U)\cong\mathbb{R}^n$ such that the numbers $\langle u,x_i\rangle$ for $1\leq i\leq r$ are distinct. Relabel the $x_i$ so that $\langle u,x_1\rangle<\langle u,x_2\rangle<\cdots<\langle u,x_r\rangle$. Each $x_i$ corresponds to some maximal cone $\sigma_i\in\mathbb{N}(\Delta)$. For $1\leq i\leq r$, define cones $s_i\in\mathcal{N}(\Delta)$ by
	\[s_i=\begin{cases}\bigcap\limits_{\substack{j>i\\ \dim(\sigma_i\cap\sigma_j)=n-1}}\sigma_i\cap\sigma_j&\text{if }1\leq i<r\\
		\sigma_r&\text{if }i=r\end{cases}.\]
The collection of cones $S_u=\{s_i\ |\ 1\leq i\leq r\}$ is called a \emph{shelling} of the fan. Inductively define sets $\Omega_i\subset X$ for $1\leq i\leq r$ by
\begin{align*}
\Omega_1&=\mathcal{O}_{s_r},\\
\Omega_i&=\bigcup_{s_{r-i+1}\subset\gamma\in\mathcal{N}(\Delta)}\mathcal{O}_\gamma\setminus\Omega_{i-1}\text{ when }1<i\leq r.
\end{align*}
The collection of $\Omega_i$'s form a CW decomposition of $X$, and $\dim\Omega_i\leq\dim\Omega_j$ whenever $i\leq j$. Hence, to a generic vector $u\in\mathbb{R}^n$, we can associate a CW decomposition of $X$. This association is locally constant in $u$. Indeed, the normal hyperplanes to the vectors $x_i-x_j$ for $1\leq i<j\leq r$ divide $\mathbb{R}^n$ into chambers, and the shelling depends only on which chamber $u$ belongs to. Equipping $X$ with a CW decomposition allows us to define various bases for $\mathrm{H}^*(X)$ and $\mathrm{K}^0(X)$.

\subsection{Bases of Characteristic Classes}
It is easy to see from the orbit structure of $X$ that $\forall 1\leq i\leq r\ \ \overline{\Omega}_i=V_{s_i}$. Thus, the cohomological and K-theoretic fundamental classes $[V_{s_i}]$ form a basis for cohomology and K-theory respectively. It follows from the triangularity of $\mathrm{c^{sm}}(\Omega_i)$ and $\mathrm{ch}(\mathrm{mC}(\Omega_i))$ with respect to cohomological grading that the $h$-deformed classes also form a basis. Hence, we may ask for the structure constants associated to any one of these bases. The structure constants associated with the basis of cohomological fundamental classes are well understood in the context of intersection theory. See \cite{Fu2} for general results in this direction. When $\mathcal{N}(\Delta)$ is generated by the Weyl chambers of a root system, a combinatorial rule for the nonequivariant structure constants associated with the basis of cohomological fundamental classes is given in \cite{Abe}. Simple examples of cohomological $h$-deformed structure constants are computed in the next two sections.

\subsection{Example: The Hirzebruch Surface}
The Hirzebruch surface $\mathcal{H}$ is the blowup of $\mathbb{P}^2$ at a $T$-fixed point. Its moment polytope $\Delta$ is the trapezoid with vertices $(0,0),(1/2,0),(0,1),(1/2,1/2)$. Let $v_1=e_1,v_2=e_2,v_3=-e_1,v_4=-e_1-e_2\in\mathbb{R}^2$. Let $\tau_i=\mathbb{R}_{\geq 0}\cdot v_i$ for $i=1,2,3,4$. Recall that given $A\subset[1,4]$, we define $\gamma_A=\mathbb{R}_{\geq 0}\cdot\{v_i\ | \ i\in A\}$. Then,
	\[\mathcal{N}(\Delta)=\{\{0\},\tau_1,\tau_2,\tau_3,\tau_4,\sigma_1:=\gamma_{\{1,4\}},\sigma_2:=\gamma_{\{1,2\}},\sigma_3:=\gamma_{\{3,4\}},\sigma_4:=\gamma_{\{2,3\}}\}.\]
The novel feature of $\mathcal{H}$ is that not all of its toric divisors are equivalent. In fact, while $[V_{\tau_2}]=[V_{\tau_4}]\ \ (*)$, we have the relation $[V_{\tau_1}]=[V_{\tau_3}]+[V_{\tau_4}]\ \ (**)$. Let $u=2e_1-e_2$. Then, the corresponding shelling $S_u$ consists of the four cones $s_1=\{0\},s_2=\tau_2,s_3=\tau_3,s_4=\sigma_4$. This gives the following CW decomposition.
	\[\Omega_1=\mathcal{O}_{\sigma_4},\ \ \Omega_2=\mathcal{O}_{\tau_3}\cup\mathcal{O}_{\sigma_3},\ \ \Omega_3=\mathcal{O}_{\tau_2}\cup\mathcal{O}_{\sigma_2},\ \ \Omega_4=\mathcal{O}_{\{0\}}\cup\mathcal{O}_{\tau_1}\cup\mathcal{O}_{\tau_4}\cup\mathcal{O}_{\sigma_1}\]
By the additivity of $\mathrm{c^{sm}}$ classes and the formula from Section 2,
\begin{align*}
\mathrm{c^{sm}}(\Omega_1)&=[V_{\tau_2}][V_{\tau_3}],\\
\mathrm{c^{sm}}(\Omega_2)&=[V_{\tau_3}][V_{\tau_4}]+h[V_{\tau_3}],\\
\mathrm{c^{sm}}(\Omega_3)&=[V_{\tau_1}][V_{\tau_2}]+h[V_{\tau_2}],\\
\mathrm{c^{sm}}(\Omega_4)&=[V_{\tau_1}][V_{\tau_4}]+h[V_{\tau_1}]+h[V_{\tau_4}]+h^2.
\end{align*}
The structure constants can be computed using relations $(*)$ and $(**)$. They are listed in Table 1.

\begin{table}[h!]
\caption{Cohomological structure constants of the Hirzebruch surface.}
\begin{equation*}
\begin{array}{ll}
\mathrm{c^{sm}}(\Omega_1)\mathrm{c^{sm}}(\Omega_1)=0&\mathrm{c^{sm}}(\Omega_2)\mathrm{c^{sm}}(\Omega_2)=-h^2\mathrm{c^{sm}}(\Omega_1)\\

\mathrm{c^{sm}}(\Omega_1)\mathrm{c^{sm}}(\Omega_2)=0&\mathrm{c^{sm}}(\Omega_2)\mathrm{c^{sm}}(\Omega_3)=h^2\mathrm{c^{sm}}(\Omega_1)\\

\mathrm{c^{sm}}(\Omega_1)\mathrm{c^{sm}}(\Omega_3)=0&\mathrm{c^{sm}}(\Omega_2)\mathrm{c^{sm}}(\Omega_4)=h^2(\mathrm{c^{sm}}(\Omega_2)+\mathrm{c^{sm}}(\Omega_1))\\

\mathrm{c^{sm}}(\Omega_1)\mathrm{c^{sm}}(\Omega_4)=h^2
\end{array}
\end{equation*}
\begin{equation*}
\begin{array}{l}
\mathrm{c^{sm}}(\Omega_3)\mathrm{c^{sm}}(\Omega_3)=0\\

\mathrm{c^{sm}}(\Omega_3)\mathrm{c^{sm}}(\Omega_4)=h^2(\mathrm{c^{sm}}(\Omega_3)+\mathrm{c^{sm}}(\Omega_1))\\

\mathrm{c^{sm}}(\Omega_4)\mathrm{c^{sm}}(\Omega_4)=h^2(\mathrm{c^{sm}}(\Omega_4)+2\mathrm{c^{sm}}(\Omega_3)+\mathrm{c^{sm}}(\Omega_2)+\mathrm{c^{sm}}(\Omega_1))
\end{array}
\end{equation*}
\end{table}
Notice that the structure constants are integer multiples of $h^2=h^{\dim(\mathcal{H})}$, one of which is negative. Moreover, the collection of structure constants appears to be independent of the choice of shelling.

\subsection{Example: The $A_2$ Permutohedral Variety}
Following \cite{Abe}, we consider the toric variety $X$ whose normal fan is given by Weyl chambers of the $A_2$ root system. The usual conventions for toric manifolds and $A_2$ are incompatible, as one of the fundamental coweights will be a vector with irrational slope. This can be remedied by defining the toric manifold with respect to the coweight lattice rather than the standard lattice. Let $E=\{(x,y,z)\in\mathbb{R}^3\ |\ x+y+z=0\}$. Define the root system $A_2=\{e_i-e_j\ |\ 1\leq i,j\leq 3\}$ on $E$. Let $\alpha_1=e_1-e_2,\alpha_2=e_2-e_3$. Then, $\Pi=\{\alpha_1,\alpha_2\}$ is a set of simple roots. The Weyl group $W=S_3$ acts on $\Phi$ by permuting indices. That is, $\forall w\in W\ \ w(e_i-e_j)=e_{w(i)}-e_{w(j)}$. The fundamental coweights are $\omega_1=2/3e_1-1/3e_2-1/3e_3,\omega_2=1/3e_1+1/3e_2-2/3e_3$, where we have identified $E$ with $E^*$ via the standard inner product on $\mathbb{R}^3$ restricted to $E$. Define a linear isomorphism $p:E\to\mathbb{R}^2$ by $p(\omega_1)=e_1,p(\omega_2)=e_2$. Abusing notation, we will denote the standard basis vectors of $\mathbb{R}^2$ by $\omega_1,\omega_2$. Taking the image of Weyl chambers and their proper faces under $p$ yields a fan of cones in $\mathbb{R}^2$ that is smooth (corresponding to a Delzant polytope). We can index the Weyl chambers by elements $w\in W$ by
	\[C_w=\mathbb{R}_{\geq 0}\cdot\{w(\omega_1),w(\omega_2)\}.\]
Let $\sigma_w=p(C_w)$. The maximal cones of the fan are
\begin{equation*}
\begin{array}{lcllcl}
\sigma_{\mathrm{id}}&=&\mathbb{R}_{\geq 0}\cdot\{\omega_1,\omega_2\}&\sigma_{(1,2)}&=&\mathbb{R}_{\geq 0}\cdot\{\omega_2-\omega_1,\omega_2\}\\
\sigma_{(1,3)}&=&\mathbb{R}_{\geq 0}\cdot\{-\omega_1,-\omega_2\}&\sigma_{(2,3)}&=&\mathbb{R}_{\geq 0}\cdot\{\omega_1,\omega_1-\omega_2\}\\
\sigma_{(1,2,3)}&=&\mathbb{R}_{\geq 0}\cdot\{-\omega_1,\omega_2-\omega_1\}&\sigma_{(1,3,2)}&=&\mathbb{R}_{\geq 0}\cdot\{-\omega_2,\omega_1-\omega_2\}
\end{array}.
\end{equation*}
The moment polytope $\Delta$ is a hexagon with vertices $(0,0),(3/4,0),(1,1/4),(1,1),(1/4,1),(0,3/4)$. This is obtained from a square by cutting off the upper left and bottom right corners. In particular, $X$ is the blowup of $\mathbb{P}^1\times\mathbb{P}^1$ at two fixed points. The fan $\mathcal{N}(\Delta)$ contains the $\sigma_w$ for all $w\in W$ along with their proper faces. For $w\in W$ and $i=1,2$, let $\tau_{w,i}=p(w(\omega_i))$. The rays in $\mathcal{N}(\Delta)$ are
	\[\mathcal{N}(\Delta)=\{\{0\}\}\cup\{\sigma_w\ |\ w\in W\}\cup\{\tau_{w,i}\ |\ w\in W,i=1,2\}.\]

Fix $u=-\omega_1-3\omega_2$. This gives an ordering of the vertices of $\Delta$ that induces the ordering
	\[\mathrm{id}<(1,2)<(1,2,3)<(2,3)<(1,3,2)<(1,3)\]
on $W$. The shelling $S_u$ consists of these six cones, which we may also index by $W$.
\begin{align*}
&s_{\mathrm{id}}:=s_1=\{0\}\ ,\ s_{(1,2)}:=s_2=\tau_{(1,2),1}\ ,\ s_{(1,2,3)}:=s_3=\tau_{(1,2,3),2}\ ,\\
&s_{(2,3)}:=s_4=\tau_{(2,3),2}\ ,\ s_{(1,3,2)}:=s_5=\tau_{(1,3,2),1}\ ,\ s_{(1,3)}:=s_6=\sigma_{(1,3)}
\end{align*}
We have chosen this particular shelling, because it conforms with the descents of the Weyl group elements. Our choice is consistent with Section 4 of \cite{Abe}. The corresponding cells can also be indexed by elements of $W$. The cells are
\begin{equation*}
\begin{array}{lcllcl}
\Omega_{\mathrm{id}}&=&\mathcal{O}_{\{0\}}\cup\mathcal{O}_{\tau_{\mathrm{id},1}}\cup\mathcal{O}_{\tau_{\mathrm{id},2}}\cup\mathcal{O}_{\sigma_{\mathrm{id}}}&\Omega_{(1,2)}&=&\mathcal{O}_{\tau_{(1,2),1}}\cup\mathcal{O}_{\sigma_{(1,2)}}\\
\Omega_{(1,3)}&=&\mathcal{O}_{\sigma_{(1,3)}}&\Omega_{(2,3)}&=&\mathcal{O}_{\tau_{(2,3),2}}\cup\mathcal{O}_{\sigma_{(2,3)}}\\
\Omega_{(1,2,3)}&=&\mathcal{O}_{\tau_{(1,2,3),2}}\cup\mathcal{O}_{\sigma_{(1,2,3)}}&\Omega_{(1,3,2)}&=&\mathcal{O}_{\tau_{(1,3,2),1}}\cup\mathcal{O}_{\sigma_{(1,3,2)}}
\end{array}.
\end{equation*}
The procedure for calculating the $\mathrm{c^{sm}}$ classes of the cells in the associated CW decomposition is the same as that of the previous example.
A few representative structure constants are listed in Table 2.

\begin{table}[h!]
\caption{Cohomological structure constants of the $A_2$ toric surface.}
\begin{align*}
\mathrm{c^{sm}}(\Omega_{(1,3,2)})\mathrm{c^{sm}}(\Omega_{(1,3,2)})&=-h^2\mathrm{c^{sm}}(\Omega_{(1,3)})\\
\mathrm{c^{sm}}(\Omega_{(1,3,2)})\mathrm{c^{sm}}(\Omega_{(2,3)})&=h^2\mathrm{c^{sm}}(\Omega_{(1,3)})\\
\mathrm{c^{sm}}(\Omega_{\mathrm{id}})\mathrm{c^{sm}}(\Omega_{(1,2)})&=h^2(\mathrm{c^{sm}}(\Omega_{(1,2)})+\mathrm{c^{sm}}(\Omega_{(1,3)}))\\
\mathrm{c^{sm}}(\Omega_{\mathrm{id}})\mathrm{c^{sm}}(\Omega_{(1,2,3)})&=h^2\mathrm{c^{sm}}(\Omega_{(1,2,3)})\\
\mathrm{c^{sm}}(\Omega_{\mathrm{id}})\mathrm{c^{sm}}(\Omega_{\mathrm{id}})&=h^2(\mathrm{c^{sm}}(\Omega_{\mathrm{id}})-\mathrm{c^{sm}}(\Omega_{(1,2,3)})-\mathrm{c^{sm}}(\Omega_{(1,3,2)})+3\mathrm{c^{sm}}(\Omega_{(1,3)}))
\end{align*}
\end{table}

\section{Grassmannians}
\subsection{Weight Function Orthogonality}
Fix $d\leq n\in\mathbb{N}$. Let $T=(\mathbb{C}^*)^n$ act on $\mathbb{C}^n$ in the usual way. This induces an action of $T$ on $\mathrm{Gr}(d,n)$. Let $e_1,...,e_n$ be the standard basis for $\mathbb{C}^n$. Then, the fixed points of $\mathrm{Gr}(d,n)$ under the $T$ action are those subspaces spanned by a subset of $\{e_1,...,e_n\}$ of cardinality $d$. The $B$-orbits of the $\binom{n}{d}$ fixed points give a decomposition of the Grassmannian into Schubert cells. Let $\mathcal{I}=\{I\subset[1,n]\ |\ |I|=d\}$. Subsets $I=\{i_1,...,i_d\}\in\mathcal{I}$ parametrize both the fixed points and the Schubert cells of $\mathrm{Gr}(d,n)$, as described above. By convention, we will enumerate the elements of $I$ in ascending order, i.e $i_1<i_2<\cdots<i_d$. Denote the fixed point $\mathrm{span}\{e_{i_1},...,e_{i_d}\}\in\mathrm{Gr}(d,n)$ by $x_I$ and the corresponding Schubert cell by $\Omega_I$. Define a partial order on $\mathcal{I}$ by $I\leq J \Leftrightarrow\forall (1\leq a\leq d)\ \  i_a\leq j_a$. This gives the Bruhat order on the Schubert cells. Another description of fixed points and Schubert cells is also commonly used.

There is a correspondence between $\mathcal{I}$ and integer partitions with $d$ parts each at most $n-d$. Such partitions are represented by Young diagrams contained in a $d\times(n-d)$ box. The correspondence is given by
	\[I=\{i_1,...,i_d\}\mapsto \lambda_I=((n-d)-(i_a-a))_{a=1}^d.\]
Under this correspondence, the number of boxes in the Young diagram of $\lambda_I$ is the codimension of the Schubert cell $\Omega_I$. This description is included for reference, but will not be used in what follows.

We will use weight function orthogonality \cite{RTV1,RTV2,TV} to obtain formulas for equivariant structure constants. See \cite{R} for a survey of generalized structure constants of flag manifolds from the perspective of weight function orthogonality. The nonequivariant structure constants can then be obtained by taking a limit using residue calculus. This residue calculus approach is also used in \cite{FR} to obtain formulas for the CSM classes of degeneracy loci. To this end, it is not necessary to fully describe the equivariant cohomology ring $\mathrm{H}^*_T(\mathrm{Gr}(d,n))$. Instead, we simply remark that the inclusion $\mathrm{Gr}(d,n)^T\hookrightarrow\mathrm{Gr}(k,n)$ induces an injection 
	\[\mathrm{H}^*_T(\mathrm{Gr}(d,n))\hookrightarrow\mathrm{H}^*_T(\mathrm{Gr}(d,n)^T)\cong\bigoplus_{I\in\mathcal{I}}\mathrm{H}^*_T(x_I)\cong\bigoplus_{I\in\mathcal{I}}\mathbb{Z}[z_1,...,z_n].\]
It follows that an equivariant cohomology class is identified with a tuple $(f_I)_{I\in\mathcal{I}}$ of polynomials in $n$ variables, one polynomial for each fixed point. The same is true in K-theory, except polynomials are replaced by Laurent polynomials. Namely,
	\[\mathrm{K}^0_T(\mathrm{Gr}(d,n))\hookrightarrow\mathrm{K}^0_T(\mathrm{Gr}(d,n)^T)\cong\bigoplus_{I\in\mathcal{I}}\mathrm{K}^0_T(x_I)\cong\bigoplus_{I\in\mathcal{I}}\mathbb{Z}[z_1^{\pm 1},...,z_n^{\pm 1}].\]
Weight functions compute the fixed point restrictions of the equivariant $\mathrm{c^{sm}}$ and $\mathrm{mC}$ classes. The key formulas are given below. Note that while weight functions are formulated for arbitrary partial flag manifolds, the formulas below have been specialized to Grassmannians.

\begin{weight}\label{weight}
Let $a\in[1,d]$, $b\in[1,n]$. Define the following polynomials in variables $x,h$.
	\[\Psi^{\mathrm{H}}_{I,a,b}(x)=\begin{cases}x+h&\text{if }b<i_a\\h&\text{if }b=i_a\\x&\text{if }b>i_a\end{cases}\quad,\quad\Psi^{\mathrm{K}}_{I,a,b}(x)=\begin{cases}1+hx&\text{if }b<i_a\\(1+h)x&\text{if }b=i_a\\1-x&\text{if }b>i_a\end{cases}.\]
Let $I=\{i_1,...,i_d\}\in\mathcal{I}$, and $\sigma\in S_n$. Define the following rational functions in variables $t_1,...,t_d,z_1,...,z_n,h$.
\begin{align*}
U^{\mathrm{H}}_I&=\prod_{a=1}^d\prod_{b=1}^n\Psi^{\mathrm{H}}_{I,a,b}(z_b-t_a)\prod_{a<b\leq d}\frac{1}{t_b-t_a}\prod_{b\leq a\leq d}\frac{1}{t_b-t_a+h},\\
U^{\mathrm{K}}_I&=\prod_{a=1}^d\prod_{b=1}^n\Psi^{\mathrm{K}}_{I,a,b}(t_a/z_b)\prod_{a<b\leq d}\frac{1}{1-t_a/t_b}\prod_{b\leq a\leq d}\frac{1}{1+ht_a/t_b},\\
W_{\sigma,I}^{\mathrm{H}}&=\mathrm{Sym}_{t_1,...,t_d}U^{\mathrm{H}}_{\sigma^{-1}(I)}(t_1,...,t_d;z_{\sigma(1)},...,z_{\sigma(n)};h),\\
W_{\sigma,I}^{\mathrm{K}}&=\mathrm{Sym}_{t_1,...,t_d}U^{\mathrm{K}}_{\sigma^{-1}(I)}(t_1,...,t_d;z_{\sigma(1)},...,z_{\sigma(n)};h),
\end{align*}
\begin{equation*}
\begin{array}{lcl}
\begin{displaystyle}R^{\mathrm{H}}_I=\prod_{\substack{a\in I\\b\notin I}}(z_b-z_a)\end{displaystyle}&,&\begin{displaystyle}Q^{\mathrm{H}}_I=\prod_{\substack{a\in I\\b\notin I}}(z_b-z_a+h)\end{displaystyle},\\
\begin{displaystyle}R^{\mathrm{K}}_I=\prod_{\substack{a\in I\\b\notin I}}(1-z_a/z_b)\end{displaystyle}&,&\begin{displaystyle}Q^{\mathrm{K}}_I=\prod_{\substack{a\in I\\b\notin I}}(1+z_b/(hz_a)).\end{displaystyle}
\end{array}
\end{equation*}

Define inner products on $\mathbb{Z}(t_1,...,t_d;z_1,...,z_n;h)$
\begin{align*}
\langle f,g\rangle^{\mathrm{H}}&=\sum_{I\in\mathcal{I}}\frac{f(z_{i_1},...,z_{i_d};z_1,...,z_n;h)g(z_{i_1},...,z_{i_k};z_1,...,z_n;h)}{R^{\mathrm{H}}_IQ^{\mathrm{H}}_I},\\
\langle f,g\rangle^{\mathrm{K}}&=\sum_{I\in\mathcal{I}}\frac{f(z_{i_1},...,z_{i_d};z_1,...,z_n;h)g(z_{i_1},...,z_{i_k};z_1,...,z_n;h)}{R^{\mathrm{K}}_IQ^{\mathrm{K}}_I}.
\end{align*}
\end{weight}

The next theorem summarizes some of the main results of \cite{RTV1,RTV2}.

\begin{weight_orth}\label{weight_orth}
Let $s_0\in S_n$ be the longest permutation. Let $\iota:\mathbb{Z}(t_1,...,t_d;z_1,...,z_n;h)\to\mathbb{Z}(t_1,...,t_d;z_1,...,z_n;h)$ be defined by $f(t_1,...,t_d;z_1,...,z_n;h)\mapsto f(1/t_1,...,1/t_d;1/z_1,...,1/z_n;1/h)$. For all $I,J\in\mathcal{I}$,
\begin{enumerate}
\item $\mathrm{c^{sm}}(\Omega_I)|_{x_J}=W^{\mathrm{H}}_{\mathrm{id},I}(z_{j_1},...,z_{j_d};z_1,...,z_n;h)\ \  ,\ \ \mathrm{mC}(\Omega_I)|_{x_J}=W^{\mathrm{K}}_{\mathrm{id},I}(z_{j_1},...,z_{j_d};z_1,...,z_n;h),$
\item $\langle W^{\mathrm{H}}_{\mathrm{id},I},W^{\mathrm{H}}_{s_0,J}\rangle^{\mathrm{H}}=\delta_{I,J}\ \ ,\ \ \langle W^{\mathrm{K}}_{\mathrm{id},I},(-h)^{-\dim \Omega_J}\iota(W^{\mathrm{K}}_{s_0,J})\rangle^{\mathrm{K}}=\delta_{I,J}.$
\end{enumerate}
\end{weight_orth}

\subsection{Cohomological Structure Constants of $\mathrm{Gr}(d,n)$}
The classes $\mathrm{c^{sm}}(\Omega_I)$ for $I\in\mathcal{I}$ form a basis for $\mathrm{H}^*_T(\mathrm{Gr}(d,n))(h)$. Let $\hat{c}^K_{I,J}\in\mathbb{Z}[z_1,...,z_n](h)$ be the unique rational functions satisfying $\forall I,J\in\mathcal{I}\ \ \mathrm{c}_T^{\mathrm{sm}}(\Omega_I)\mathrm{c}_T^{\mathrm{sm}}(\Omega_J)=\sum_{K\in\mathcal{I}}\hat{c}^K_{I,J}\mathrm{c}_T^{\mathrm{sm}}(\Omega_K)$. The following is an immediate consequence of the orthogonality relations in Theorem~\ref{weight_orth}.

\begin{H_weight_LR}\label{H_weight_LR}
For all $I,J,K\in\mathcal{I}$,\ \ $\hat{c}^K_{I,J}=\langle W^{\mathrm{H}}_{\mathrm{id},I}W^{\mathrm{H}}_{\mathrm{id},J},W^{\mathrm{H}}_{\mathrm{s_0},K}\rangle^{\mathrm{H}}$.
\end{H_weight_LR}

A priori, this formula for the structure constants gives rational functions in $z_1,...,z_n$. They are in fact polynomials, but the cancellations of the denominators are not obvious. If we let $c^K_{I,J}=\lim_{z_1,...,z_n\to 0}\hat{c}^K_{I,J}$, then the $c^K_{I,J}$ are the structure constants in nonequivariant cohomology. It will be convenient to consider each term of the summation in Corollary~\ref{H_weight_LR} individually, so define
	\[\hat{c}^{K,L}_{I,J}=\frac{ W^{\mathrm{H}}_{\mathrm{id},I}W^{\mathrm{H}}_{\mathrm{id},J}W^{\mathrm{H}}_{\mathrm{s_0},K}(z_{i_1},...,z_{i_d};z_1,...,z_n;h)}{R^{\mathrm{H}}_IQ^{\mathrm{H}}_I}.\]
It is clear from the formulas that $\hat{c}^{K,L}_{I,J}=0$ unless $K\leq L\leq I,J$. Let us revisit $\mathbb{P}^m$ from this perspective.

\subsubsection{$\mathbb{P}^m$ Revisited}
Here, we will abuse notation by letting an integer $i\in\mathbb{Z}$ also represent the set $\{i\}$. With this notation, the equivariant structure constants of Gr$(1,n)$ are denoted $\hat{c}^k_{i,j}$. This convention is consistent with the convention of Section 2. Without loss of generality, assume $i\leq j$. Then, Corollary~\ref{H_weight_LR} gives us the following formula for $\hat{c}^{k,l}_{i,j}$, after making obvious cancellations:
	\[\hat{c}^{k,l}_{i,j}=\begin{cases}\frac{h\prod\limits_{b=j+1}^{n}(z_b-z_l)\prod\limits_{b=1}^{i-1}(z_b-z_l+h)\prod\limits_{b=k+1}^{j-1}(z_b-z_l+h)}{\prod\limits_{b=k}^{l-1}(z_b-z_l)\prod\limits_{b=l+1}^i(z_b-z_l)}&\text{if }k\leq l\leq i,\\
	0 &\text{ otherwise.}\end{cases}\]
Note that if $k=l=i$, then this formula reduces to the polynomial $h\prod\limits_{b=j+1}^{n}(z_b-z_i)\prod\limits_{\substack{b=1\\b\neq i}}^{j-1}(z_b-z_i+h)$. Hence, $c^i_{i,j}=\begin{cases}h^{n-1}&\text{if }j=n\\0&\text{otherwise}\end{cases}$. It remains to compute $c^k_{i,j}$ for $k<i$. 

Assume that $k<i$. Since $\hat{c}^k_{i,j}\in\mathrm{H}^*_T(\mathrm{Gr}(1,n))$ is a polynomial in $z_b$ variables, we may apply Cauchy's integral theorem to compute $c^k_{i,j}$. For $1\leq b\leq n$, let $\gamma_b\subset\mathbb{C}$ be a counterclockwise circle centered on 0 with radius $1/(b+1)$. Then,
\begin{align*}
c^k_{i,j}&=\left(\frac{1}{2\pi\sqrt{-1}}\right)^n\int_{\gamma_n}\cdots\int_{\gamma_1}\frac{\hat{c}^k_{i,j}}{z_1\cdots z_n}dz_1\cdots dz_n\\
	&=\left(\frac{1}{2\pi\sqrt{-1}}\right)^n\int_{\gamma_n}\cdots\int_{\gamma_1}\sum_{l=1}^n\frac{\hat{c}^{k,l}_{i,j}}{z_1\cdots z_n}dz_1\cdots dz_n\\
	&=\left(\frac{1}{2\pi\sqrt{-1}}\right)^n\sum_{l=1}^n\int_{\gamma_n}\cdots\int_{\gamma_1}\frac{\hat{c}^{k,l}_{i,j}}{z_1\cdots z_n}dz_1\cdots dz_n.
\end{align*}
Define $c^{k,l}_{i,j}=\left(\frac{1}{2\pi\sqrt{-1}}\right)^n\int_{\gamma_n}\cdots\int_{\gamma_1}\frac{\hat{c}^{k,l}_{i,j}}{z_1\cdots z_n}dz_1\cdots dz_n$. Our goal is to evaluate this integral using residue calculus. By Fubini's theorem, the order of integration can be freely interchanged. Define the sets of variables
	\[\mathbf{z}_1=\{z_b\ |\ b\notin[k,i]\}\ ,\ \mathbf{z}_2=\{z_b\ |\ b\in[l+1,i]\}\ ,\ \mathbf{z}_3=\{z_b\ |\ b\in[k,l-1]\}.\]
The most convenient order of integration is first over variables in $\mathbf{z}_1$, second over variables in $\mathbf{z}_2$, third over variables in $\mathbf{z}_3$, and finally over $z_l$. From now on, we do not fully notate this iterated integral. Instead, the symbol $\int fd\mathbf{z}_s$, where $s=1,2,3$, will represent the iterated integral of $\left(\frac{1}{2\pi\sqrt{-1}}\right)f$ taken over the variables in $\mathbf{z}_s$.

The first set of integrals $\int\frac{\hat{c}^{k,l}_{i,j}}{z_1\cdots z_n}d\mathbf{z}_1$ picks up singularities only at $z_b=0$ for $b\notin[k,i]$. Since $\hat{c}^{k,l}_{i,j}$ is holomorphic in $z_b$ when $b\notin[k,i]$, this integral just sets these variables equal to 0. We thus have
	\[\int\frac{\hat{c}^{k,l}_{i,j}}{z_1\cdots z_n}d\mathbf{z}_1=-\frac{h(-z_l+h)^{k-1+(j-i-1)(1-\delta_{i,j})}(z_k-z_l+h)(z_i-z_l+h)^{1-\delta_{i,j}}\prod\limits_{b=k+1}^{i-1}(z_b-z_l+h)^2}{(-z_l)^{j-n+1}\prod\limits_{b=k}^{l-1}(z_b-z_l)z_b\prod\limits_{b=l+1}^i(z_b-z_l)z_b}.\]
Due to the way the $\gamma_b$ are nested, the second set of integrals $\int\int\frac{\hat{c}^{k,l}_{i,j}}{z_1\cdots z_n}d\mathbf{z}_1 d\mathbf{z}_2$ picks up singularities only at $z_b=0$ for $b\in[l+1,i]$. Similarly, we set these variables equal to 0 and obtain
	\[\int\int\frac{\hat{c}^{k,l}_{i,j}}{z_1\cdots z_n}d\mathbf{z}_1 d\mathbf{z}_2=
	\begin{cases}-\frac{h^2(-z_l+h)^{2i-k-3+(j-i)(1-\delta_{i,j})}}{(-z_l)^{i+j-n-k+1}}&\text{if }l=k,\\
	-\frac{h^3(-z_l+h)^{2i+k-2l-3+(j-i)(1-\delta_{i,j})}(z_k-z_l+h)\prod\limits_{b=k+1}^{l-1}(z_b-z_l+h)^2}{(-z_l)^{i+j-l-n+1}\prod\limits_{b=k}^{l-1}(z_b-z_l)z_b}&\text{if }k<l<i,\\
	-\frac{h^{2-\delta_{i,j}}(-z_i+h)^{k-1+(j-i-1)(1-\delta_{i,j})}(z_k-z_i+h)\prod\limits_{b=k+1}^{i-1}(z_b-z_i+h)^2}{(-z_i)^{j-n+1}\prod\limits_{b=k}^{i-1}(z_b-z_i)z_b}&\text{if }l=i.
	\end{cases}\]
We can simplify these expressions by noticing that $(j-i)(1-\delta_{i,j})=j-i$. The following residue calculations will be useful for performing the $\mathbf{z}_3$ integral.

\begin{res}\label{res}
Let $f=\frac{(z_b-z_l+h)^r}{(z_b-z_l)z_b}$, where $r=1,2$. Then,
\begin{enumerate}
\item $\mathrm{Res}_{z_b= 0}f=\frac{(-z_l+h)^r}{-z_l}$,
\item $\mathrm{Res}_{z_b= z_l}f=\frac{h^r}{z_l}$,
\item $\mathrm{Res}_{z_b= 0}f+\mathrm{Res}_{z_b= z_l}f=\begin{cases}1&\text{if }r=1,\\(-z_l+2h)&\text{if }r=2.\end{cases}$
\end{enumerate}
\end{res}
It immediately follows that
	\[\int\int\int\frac{\hat{c}^{k,l}_{i,j}}{z_1\cdots z_n}d\mathbf{z}_1 d\mathbf{z}_2 d\mathbf{z}_3=
	\begin{cases}-\frac{h^2(-z_k+h)^{i+j-k-3}}{(-z_l)^{i+j-k-n+1}}&\text{if }l=k,\\
	-\frac{h^3(-z_k+h)^{i+j+k-2l-3}(-z_l+2h)^{l-k-1}}{(-z_l)^{i+j-l-n+1}}&\text{if }k<l<i,\\
	-\frac{h^{2-\delta_{i,j}}(-z_k+h)^{k-1+(j-i-1)(1-\delta_{i,j})}(-z_l+2h)^{i-k-1}}{(-z_l)^{j-n+1}}&\text{if }l=i.
	\end{cases}\]
Finally, we perform the change of variables $z=-z_l$ and integrate over the new variable $z$. The following proposition summarizes the calculations of this section.

\begin{csm_proj_weight}\label{csm_proj_weight}
For all $i,j,k,l\in[1,n]$,
\begin{enumerate}
\item $c^{k,l}_{i,j}=0$ unless $k\leq l\leq i,j$,
\item when $k<i$ and $k\leq l\leq i\leq j$,
	\[c^{k,l}_{i,j}=
	\begin{cases}\mathrm{Res}_{z= 0}\left(\frac{h^2(z+h)^{i+j-k-3}}{z^{i+j-k-n+1}}\right)&\text{if }l=k,\\
	\mathrm{Res}_{z= 0}\left(\frac{h^3(z+h)^{i+j+k-2l-3}(z+2h)^{l-k-1}}{z^{i+j-l-n+1}}\right)&\text{if }k<l<i,\\
	\mathrm{Res}_{z= 0}\left(\frac{h^{2-\delta_{i,j}}(z+h)^{k-1+(j-i-1)(1-\delta_{i,j})}(z+2h)^{i-k-1}}{z^{j-n+1}}\right)&\text{if }l=i,
	\end{cases}\]
\item when $k=l=i\leq j$, $c^{i,i}_{i,j}=\begin{cases}h^{n-1}&\text{if }j=n,\\0&\text{otherwise}.\end{cases}$
\end{enumerate}
\end{csm_proj_weight}
Observe that each of these residues is a nonnegative integer multiple of $h^{n-1}=h^{\dim(\mathrm{Gr}(1,n))}$. Moreover, we can realize $c^k_{i,j}$ as a coefficient of some polynomial. Namely, for $1\leq k\leq l\leq i\leq j\leq n$ and $k<i$, define the polynomial
	\begin{multline*}p^k_{i,j}(z)=h^2(z+h)^{i+j-k-3}+h^{2-\delta_{i,j}}z^{i-k}(z+h)^{k-1+(j-i-1)(1-\delta_{i,j})}(z+2h)^{i-k-1}\\
		+\sum_{l=k+1}^{i-1}h^3z^{l-k}(z+h)^{i+j+k-2l-3}(z+2h)^{l-k-1}.\end{multline*}
Then $c^k_{i,j}$ is the coefficient in $p^k_{i,j}$ of the degree $i+j-k-n$ term in $z$. Applying Theorem~\ref{HLR_proj} yields nonobvious properties of these polynomials. It would be interesting to see more general results computing structure constants of Grassmannians as coefficients of explicit polynomials.

\subsubsection{Pieri Triples in Gr$(2,n)$}
The results of the previous section can be generalized to a special class of structure constants in 2-plane Grassmannians. The novel feature of 2-plane Grassmannians is the appearence of symmetrizations in the weight functions. This section focuses on the situation when these symmetrizations contribute only one nonzero term. Assume that $I,J,K\in\mathcal{I}$ are such that $n\in I,J,K$. Note that $\hat{c}^{K,L}_{I,J}=0$ unless $n\in L$ as well. Subject to this condition, $I,J,K,L$ are determined by $i:=i_1,j:=j_1,k:=k_1,l:=l_1$. Such $\Omega_I,\Omega_J,\Omega_K$ are Pieri cells, so we expect $\hat{c}^K_{I,J}$ to be related to the structure constant $\hat{c}^k_{i,j}$ of Gr$(1,n-1)=\mathbb{P}^{n-2}$. See for instance the recursions of \cite{AM2009}. Without loss of generality, assume that $i\leq j$. By making obvious cancellations in Corollary~\ref{H_weight_LR}, we get
	\[\hat{c}^{K,L}_{I,J}=\begin{cases}\left[\frac{h\prod\limits_{b=j+1}^{n-1}(z_b-z_l)\prod\limits_{b=1}^{i-1}(z_b-z_l+h)\prod\limits_{b=k+1}^{j-1}(z_b-z_l+h)}{\prod\limits_{b=k}^{l-1}(z_b-z_l)\prod\limits_{b=l+1}^i(z_b-z_l)}\right]\prod\limits_{\substack{1\leq b\leq n-1\\b\neq l}}(z_b-z_n+h) &\text{ if }k<j,\\
	\left[\prod\limits_{b=1}^{j-1}(z_b-z_l+h)\prod\limits_{b=j+1}^{n-1}(z_b-z_l)\right]\prod\limits_{\substack{1\leq b\leq n-1\\b\neq l}}(z_b-z_n+h) &\text{ if }k=l=i=j,\\
	0 &\text{ otherwise}.\end{cases}\]
The terms in square brackets are precisely the $\hat{c}^{k,l}_{i,j}$, so we indeed have a relationship between the structure constants of Gr$(2,n)$ and those of Gr$(1,n-1)$.

\begin{proj_like}\label{proj_like}
If $n\in I,J,K,L$, and $i=i_1,j=j_1,k=k_1,l=l_1$ are such that $i\leq j$, then
	\[\hat{c}^{K,L}_{I,J}=\hat{c}^{k,l}_{i,j}\prod\limits_{\substack{1\leq b\leq n-1\\b\neq l}}(z_b-z_n+h).\]
\end{proj_like}

The extra factors in Proposition~\ref{proj_like} will have a predictable effect on the residues computed in the previous section, resulting in a formula for $c^K_{I,J}$. We will adopt the notations and conventions of the previous section in what follows. In particular, define the subsets of $\{z_1,...,z_n\}$
	\[\mathbf{z}_1=\{z_b\ |\ b\notin[k,i]\}\ ,\ \mathbf{z}_2=\{z_b\ |\ b\in[l+1,i]\}\ ,\ \mathbf{z}_3=\{z_b\ |\ b\in[k,l-1]\}.\]
We will again integrate over $\mathbf{z}_1,\mathbf{z}_2,\mathbf{z}_3$, and finally $z_l$. The $d\mathbf{z}_1$ and $d\mathbf{z}_2$ integrals pick up only singularities at 0. Like before, we set variables equal to 0 to obtain
	\[\int\int\frac{\hat{c}^{K,L}_{I,J}}{z_1\cdots z_n}d\mathbf{z}_1 d\mathbf{z}_2=h^{n+k-l-2}\prod\limits_{k\leq b\leq l-1}(z_b+h)\int\int\frac{\hat{c}^{k,l}_{i,j}}{z_1\cdots z_{n-1}}d\mathbf{z}_1 d\mathbf{z}_2.\]
We have already computed $\int\int\frac{\hat{c}^{k,l}_{i,j}}{z_1\cdots z_{n-1}}d\mathbf{z}_1 d\mathbf{z}_2$ in the previous section. The following lemma will be useful for computing the $\mathbf{z}_3$ integral.

\begin{res_Pieri}\label{res_Pieri}
Let $f=(z_b+h)\frac{(z_b-z_l+h)^r}{(z_b-z_l)z_b}$, where $r=1,2$. Then,
\begin{enumerate}
\item $\mathrm{Res}_{z_b= 0}f=\frac{h(-z_l+h)^r}{-z_l}$,
\item $\mathrm{Res}_{z_b= z_l}f=\frac{(z_l+h)h^r}{z_l}$,
\item $\mathrm{Res}_{z_b= 0}f+\mathrm{Res}_{z_b= z_l}f=\begin{cases}2h&\text{if }r=1,\\h(-z_l+3h)&\text{if }r=2\end{cases}$.
\end{enumerate}
\end{res_Pieri}

It immediately follows that for $k\leq l\leq i\leq j$ and $k<i$,
	\[\int\int\int\frac{\hat{c}^{K,L}_{I,J}}{z_1\cdots z_n}d\mathbf{z}_1 d\mathbf{z}_2\mathbf{z}_3=
	\begin{cases}-\frac{h^n(-z_l+h)^{i+j-k-3}}{-z_l^{i+j-k-n+2}}&\text{if }l=k,\\
	-\frac{2h^{n+1}(-z_l+h)^{i+j+k-2l-3}(-z_l+3h)^{l-k-1}}{-z_l^{i+j-l-n+2}}&\text{if }k<l<i,\\
	-\frac{2h^{n-\delta_{i,j}}(-z_l+h)^{k-1+(j-i-1)(1-\delta_{i,j})}(-z_l+3h)^{i-k-1}}{-z_l^{j-n+2}}&\text{if }l=i.
	\end{cases}\]
Like before, we make the change of variables $z=-z_l$.

\begin{csm_Pieri}\label{csm_Pieri}
Let $I,J,K,L\in\mathcal{I}$ be such that $n\in I,J,K$. Let $i=i_1,j=j_1,k=k_1,l=l_1$. Then,
\begin{enumerate}
\item $c^{K,L}_{I,J}=0$ unless $n\in L$ and $k\leq l\leq i,j$,
\item when $k<i$, $k\leq l\leq i\leq j$, and $n\in L$,
	\[c^{K,L}_{I,J}=
	\begin{cases}\mathrm{Res}_{z= 0}\left(\frac{h^n(z+h)^{i+j-k-3}}{z^{i+j-k-n+2}}\right)&\text{if }l=k,\\
	\mathrm{Res}_{z= 0}\left(\frac{2h^{n+1}(z+h)^{i+j+k-2l-3}(z+3h)^{l-k-1}}{z^{i+j-l-n+2}}\right)&\text{if }k<l<i,\\
	\mathrm{Res}_{z= 0}\left(\frac{2h^{n-\delta_{i,j}}(z+h)^{k-1+(j-i-1)(1-\delta_{i,j})}(z+3h)^{i-k-1}}{z^{j-n+2}}\right)&\text{if }l=i.
	\end{cases}\]
\item when $k=l=i\leq j$, $c^{I,I}_{I,J}=\begin{cases}h^{2(n-2)}&\text{if }j=n-1,\\0&\text{otherwise}.\end{cases}$
\end{enumerate}
\end{csm_Pieri}
Again, we see that the structure constants are nonnegative integer multiples of $h^{2(n-2)}=h^{\dim(\mathrm{Gr}(2,n))}$. It also possible to realize $c^K_{I,J}$, where $k<i\leq j$, as the degree $i+j-k-n+1$ term of the polynomial
	\begin{multline*}p^K_{I,J}=h^n(z+h)^{i+j-k-3}+2h^{n+1}\sum_{l=k+1}^{i-1}z^{l-k}(z+h)^{i+j+k-2l-3}(z+3h)^{l-k-1}\\
		+2h^{n-\delta_{i,j}}z^{i-k}(z+h)^{k-1+(j-i-1)(1-\delta_{i,j})}(z+3h)^{i-k-1}.\end{multline*}

More generally, consider $\mathrm{Gr}(d,n)$ for arbitrary $d$. Given $I\in\mathcal{I}$ such that $n\in I$, let $I_-=I\setminus\{n\}$. The structure constants $\hat{c}^K_{I,J}$ associated to a Pieri triple $I,J,K\in\mathcal{I}$, where $n\in I,J,K$, are related to the structure constants $\hat{c}^{K_-}_{I_-,J_-}$ of $\mathrm{Gr}(d-1,n-1)$ by formulas akin to Proposition~\ref{proj_like}. However, the influence of the extra terms on the residue calculus of $\hat{c}^{K_-}_{I_-,J_-}$ may not be as straightforward.

\subsubsection{The General Case}
In general, the formula of Corollary~\ref{H_weight_LR} is complicated by symmetrizations. Take for example, the structure constant of $\mathrm{Gr}(2,4)$
	\[\hat{c}^{\{1,2\},\{2,3\}}_{\{2,3\},\{2,3\}}=(z_4-z_2)(z_4-z_3)(h+z_1-z_2)(h+z_1-z_3)\frac{\frac{h^2(z_1-z_3)}{z_2-z_3}+\frac{h(z_1-z_2)(z_3-z_2+h)}{z_3-z_2}}{(z_1-z_2)(z_1-z_3)}.\]
In light of part 1 of Theorem~\ref{weight_orth} and the fact that fixed point restrictions are always polynomials in the $z$ variables, the singularity at $z_2=z_3$ must be removable. Indeed, this expression simpifies to
	\[\frac{h(z_1-z_2+h)^2(z_1-z_3+h)(z_4-z_2)(z_4-z_3)}{(z_1-z_2)(z_1-z_3)}.\]
Notice that the simplified expression can be written as a product of three kinds of terms:
\begin{enumerate}
\item a power of $h$,
\item terms of the form $(z_b-z_l+h)^r/(z_b-z_l)^s$, where $b\notin L$, $l\in L$, $r=0,1,2$, and $s=0,1$,
\item and terms of the form $(z_b-z_l)$, where $\forall a\in L\ \ b>a$, $l\in L$,
\end{enumerate}
where terms of type 2 and 3 with fixed $b$ and $l$ appear at most once. This is precisely the pattern we encountered in the previous sections, and it appears that the $\hat{c}^{K,L}_{I,J}$ will be sums of expressions of this form in general. The novel feature is that for all $b\notin L$, up to $d$ many poles in $z_b$ of order 1 may appear. Hence, the $z_b$ integral may contribute up to $d+1$ many residues (the extra residue comes from the singularity at $z_b=0$ introduced by Cauchy's integral theorem). Computing these integrals is a subject for future work, but certain properties can be deduced from this pattern, e.g. positivity. In the next section, we will describe conjectures involving the more general K-theoretic structure constants, which might be provable using this technique.

\subsection{K-Theoretic Structure Constants of $\mathrm{Gr}(d,n)$}
The classes $\mathrm{mC}(\Omega_I)$ for $I\in\mathcal{I}$ form a basis for $\mathrm{K}^0_T(\mathrm{Gr}(d,n))(h)$. Let $\hat{C}^K_{I,J}\in\mathbb{Z}[z^{\pm 1}_1,...,z^{\pm 1}_n](h)$ be the unique rational functions satisfying $\forall I,J\in\mathcal{I}\ \ \mathrm{mC}_T(\Omega_I)\mathrm{mC}_T(\Omega_J)=\sum_{K\in\mathcal{I}}\hat{C}^K_{I,J}\mathrm{mC}_T(\Omega_K)$. The following is an immediate consequence of the orthogonality relations in Theorem~\ref{weight_orth}.

\begin{K_weight_LR}\label{K_weight_LR}
For all $I,J,K\in\mathcal{I}$,\ \ $\hat{C}^K_{I,J}=\langle W^{\mathrm{K}}_{\mathrm{id},I}W^{\mathrm{K}}_{\mathrm{id},J},(-h)^{-\dim(\Omega_k)}\iota(W^{\mathrm{K}}_{\mathrm{s_0},K})\rangle^{\mathrm{K}}$.
\end{K_weight_LR}

As in the previous section, define
	\[\hat{C}^{K,L}_{I,J}=\frac{(-h)^{-\dim(\Omega_k)}W^{\mathrm{K}}_{\mathrm{id},I}W^{\mathrm{K}}_{\mathrm{id},J}\iota(W^{\mathrm{K}}_{\mathrm{s_0},K})(z_{i_1},...,z_{i_d};z_1,...,z_n;h)}{R^{\mathrm{K}}_IQ^{\mathrm{K}}_I}.\]
Taking the limit as $z$ variables go to 1 degenerates the equivariant structure constants to nonequivariant structure constants. In order to compare the K-theoretic structure constants to the cohomological structure constants, make the change of variables $\zeta_b=z_b-1$, for all $1\leq b\leq n$. Then, $C^K_{I,J}=\lim_{\zeta_1,...,\zeta_n\to 0}\hat{C}^K_{I,J}$. We can evaluate this limit by integrating as we did in the previous section. With the contours $\gamma_b$ of the previous section, let
	\[C^{K,L}_{I,J}=\left(\frac{1}{2\pi\sqrt{-1}}\right)^n\int_{\gamma_1}\cdots\int_{\gamma_n}\frac{\hat{C}^{K,L}_{I,J}}{\zeta_1\cdots\zeta_n}d\zeta_n\cdots d\zeta_1.\]
Then, $C^K_{I,J}=\sum_{L\in I}C^{K,L}_{I,J}$. It will also be helpful to introduce the variable $\nu=-1-h$. Let us write the K-theoretic structure constant of Gr$(2,4)$ from the previous section in these new variables. We have
	\[\hat{C}^{\{1,2\},\{2,3\}}_{\{2,3\},\{2,3\}}=\frac{\nu(\zeta_1-(\nu-1)\zeta_2+\nu)^2(\zeta_1-(\nu-1)z_3+\nu)(\zeta_4-\zeta_2)(\zeta_4-z_3)}{(\zeta_1-\zeta_2)(\zeta_1-\zeta_3)}\cdot\frac{1}{(\zeta_1+1)(\zeta_4+1)^2}.\]
Notice the similarities between $\hat{C}^{\{1,2\},\{2,3\}}_{\{2,3\},\{2,3\}}$ and $\hat{c}^{\{1,2\},\{2,3\}}_{\{2,3\},\{2,3\}}$:
\begin{enumerate}
\item instead of a power of $h$, there is a power of $\nu$,
\item instead of $z_b-z_l+h$, there is $\zeta_b-(\nu-1)\zeta_l+\nu$,
\item instead of $z_b-z_l$, there is $\zeta_b-\zeta_l$,
\item and there is an extra factor of the form $\prod_{b=1}^n (\zeta_b+1)^{r_b}$, where $r_b\in\mathbb{Z}$.
\end{enumerate}
This pattern appears to hold in general. Since the contours are circles of radius less than 1, the extra terms in 4 do not contribute any new singularities. Thus, the residue calculus of these K-theoretic expressions is not substantially more complicated than that of the cohomological expressions.

We expect general versions of the results of Section 2.3 to hold for Gr$(d,n)$. 

\begin{symmetry}
For all $I,I',J,J',K,K'\in\mathcal{I}$,
\begin{enumerate}
\item $C^K_{I,J}$ is a polynomial in $\nu$ with nonnegative coefficients,
\item $C^K_{I,J}=0$ unless $K\leq I,J$,\\
\item $C^K_{I,J}=C^K_{I',J'}$ if $i_a+j_a=i'_a+j'_a$ for all $1\leq a\leq d$,\\
\item $C^K_{I,J}=C^{K'}_{I,J'}$ if $j_a-k_a=j'_a-k'_a$ for all $1\leq a \leq d$.
\end{enumerate}
\end{symmetry}

\begin{KtoH}
The term of $C^K_{I,J}$ with lowest $\nu$-degree has coefficient $c^K_{I,J}/h^{\dim(\mathrm{Gr}(d,n))}$.
\end{KtoH}

\noindent
Department of Mathematics, University of North Carolina at Chapel Hill, USA\\
\emph{email address:} yshou@live.unc.edu

\section{References}
\begin{biblist}
\bib{Abe}{article}{
title={Young Diagrams and Intersection Numbers for Toric Manifolds associated with Weyl Chambers }, 
volume={22}, 
DOI={https://doi.org/10.37236/4307 }, 
journal={The Electronic Journal of Combinatorics}, 
author={Abe, Hiraku}, 
date={2015}
}

\bib{AM2009}{article}{
    author = {Aluffi, Paolo}
    author={Mihalcea, Leonardo Constantin},
     title = {Chern classes of Schubert cells and varieties},
   journal = {Journal of Algebraic Geometry},
    volume = {18},
      date = {2009},
     pages = {63--100},
     issn = {1056-3911},
}

\bib{AM2016}{article}{
    author = {Aluffi, Paolo},
    author={Mihalcea, Leonardo Constantin},
     title = {Chern-Schwartz-MacPherson  classes  for  Schubert  cells  in  flag  manifolds},
   journal = {Compositio Mathematica},
    volume = {152},
  number={12}
      date = {2016},
     pages = {2603--2652},
publisher={London Mathematical Society}
}

\bib{AMSS1}{article}{
    author = {Aluffi, Paolo},
    author={Mihalcea, Leonardo Constantin},
    author = {Schürmann, J.},
    author={Su, C.},
        title = {Shadows of characteristic cycles, Verma modules, and positivity of Chern-Schwartz-MacPherson classes of Schubert cells},
      journal = {arXiv e-prints},
         year = {2017},
          eid = {arXiv:1709.08697},
        pages = {arXiv:1709.08697},
archivePrefix = {arXiv},
       eprint = {1709.08697},
 primaryClass = {math.AG}
}

\bib{AMSS2}{article}{
    author = {Aluffi, Paolo},
    author={Mihalcea, Leonardo Constantin},
    author = {Schürmann, J.},
    author={Su, C.},
    title = {Motivic Chern classes of Schubert cells, Hecke algebras, and applications to Casselman's problem},
    journal = {arXiv e-prints},
    year = {2019},
          eid = {arXiv:1902.10101},
        pages = {arXiv:1902.10101},
archivePrefix = {arXiv},
       eprint = {1902.10101},
 primaryClass = {math.AG}
}

\bib{BSY}{article}{
    author = {Brasselet, Jean-Paul}
    author={Schürmann, Jörg},
    author={Yokura, Shoji},
     title = {Hirzebruch classes and motivic Chern classes for singular spaces},
   journal = {Journal of Topology and Analysis},
    volume = {2},
      date = {2010},
     pages = {1--55},
     number = {1}
}

\bib{FR}{article}{
title={Chern-Schwarz-Macpherson Classes of Degeneracy Loci}, 
volume={22}, 
journal={Geometry and Topology}, 
author={Feh\'er, László},
author={Rimányi, Richárd}, 
date={2018},
pages={3575–3622}
}

\bib{Fu}{book}{ 
title={Introduction to Toric Varieties}, 
date={1993}, 
author={Fulton, William}
publisher={Princeton University Press}
series={AM-131}
}

\bib{Fu1}{webpage}{ 
title={Equivariant Cohomology in Algebraic Geometry}, 
date={2007}, 
url={https://people.math.osu.edu/anderson.2804/eilenberg/lecture13.pdf}, 
author={Fulton, William}
}

\bib{Fu2}{article}{
title = {Intersection theory on toric varieties},
journal = {Topology},
volume = {36},
pages = {335 - 353},
year = {1997},
issn = {0040-9383},
doi = {https://doi.org/10.1016/0040-9383(96)00016-X},
author = {Fulton, William},
author = {Sturmfels, Bernd}
}

\bib{M}{article}{
author = {Macpherson, R. D.}
title = {Chern Classes of Singular Algebraic Varieties},
journal = {Annals of Mathematics},
volume = {100},
number = {2},
pages = {423-432},
year = {1974}
}

\bib{MS}{article}{
author = {Maxim, Laurenţiu G.}
author = {Schürmann, Jörg},
title = {Characteristic Classes of Singular Toric Varieties},
journal = {Communications on Pure and Applied Mathematics},
volume = {68},
number = {12},
pages = {2177-2236},
doi = {10.1002/cpa.21553},
year = {2015}
}

\bib{R}{article}{
       author = {Rimányi, Richárd},
        title = {$\hbar$-deformed Schubert calculus in equivariant cohomology, K-theory, and elliptic cohomology},
      journal = {arXiv e-prints},
     keywords = {Mathematics - Algebraic Geometry, 14N15, 55N34},
         year = {2019},
archivePrefix = {arXiv},
       eprint = {https://arxiv.org/abs/1912.13089},
 primaryClass = {math.AG},
}

\bib{RTV1}{article}{
title={Partial flag varieties, stable envelopes and weight functions}, 
volume={6}, 
journal={Quantum Topology}, 
author={Rimányi, Richárd}, 
author={Tarasov, Vitaly},
author={Varchenko, Alexander},
date={2015},
pages={333--364}
}

\bib{RTV2}{article}{
title={Trigonometric weight functions as K-theoretic stable envelope maps for the cotangent bundle of a flag variety}, 
volume={94}, 
journal={Journal of Geometry and Physics}, 
author={Rimányi, Richárd}, 
author={Tarasov, Vitaly},
author={Varchenko, Alexander},
date={2015},
pages={81--119}
}

\bib*{banach}{proceedings}{
  conference = {IMPANGA2015},
  editor =       {Buczynski, J.},
  editor =       {Michalek, M.},
  editor =       {Postingel, E.},
}

\bib{RV}{article}{
  xref={banach},
  conference = {title={IMPANGA2015}},
  author =       {Rimányi, R.},
  author =       {Varchenko, A.},
  title =        {Equivariant Chern-Schwartz-MacPherson classes in partial flag varieties: interpolation and formulae},
  book={title={in Schubert Varieties, Equivariant Cohomology and Characteristic Classes}}
  pages =        {225--235},
  date={2018}
}

\bib{S}{article}{
  author = {{Su}, Changjian},
        title = {Structure constants for Chern classes of Schubert cells},
      journal = {arXiv e-prints},
     keywords = {Mathematics - Algebraic Geometry, Mathematics - Combinatorics},
         year = {2019},
          eid = {arXiv:1909.10940},
        pages = {arXiv:1909.10940},
archivePrefix = {arXiv},
       eprint = {1909.10940},
 primaryClass = {math.AG},
}

\bib{TV}{article}{
title={Hypergeometric solutions of the quantum differential equation of the cotangent bundle of a partial flag variety}, 
volume={12}, 
journal={Central European Journal of Mathematics}, 
author={Tarasov, Vitaly},
author={Varchenko, Alexander},
date={2013},
doi = {10.2478/s11533-013-0376-8},
}
\end{biblist}
\end{document}